\documentclass[hidelinks,onefignum,onetabnum]{siamart250211}

\usepackage{lipsum}
\usepackage{amsfonts}
\usepackage{graphicx}
\usepackage{epstopdf}
\usepackage{algorithmic}
\ifpdf
  \DeclareGraphicsExtensions{.eps,.pdf,.png,.jpg}
\else
  \DeclareGraphicsExtensions{.eps}
\fi

\usepackage{amsfonts}       \usepackage{amsmath}       \usepackage{amssymb}       

\usepackage{fancyhdr}       \usepackage{graphicx}       
\usepackage{multirow}

\usepackage{pgfplots}
\usepackage{subcaption}
\usepackage{cleveref}
\usepackage{siunitx}
\usepackage{ulem}

\renewcommand{\vec}[1]{\mathbf{#1}}

\newsiamremark{remark}{Remark}
\newsiamremark{hypothesis}{Hypothesis}
\crefname{hypothesis}{Hypothesis}{Hypotheses}
\newsiamthm{claim}{Claim}
\newsiamremark{fact}{Fact}
\crefname{fact}{Fact}{Facts}

\headers{OSM vs. ORAS}{B. Martin, P. Jolivet and C. Geuzaine}

\title{Comparison of substructured non-overlapping domain decomposition and
overlapping additive Schwarz methods for large-scale Helmholtz
problems with multiple sources\thanks{Submitted on June 20th, 2025.
\funding{This work was funded by a FRIA doctoral fellowship from the Fonds de la Recherche Scientifique - FNRS, Belgium and by the Walloon Region through the Win2Wal EXPANSION project (Grant No. 2010161).}}}

\author{Dianne Doe\thanks{Imagination Corp., Chicago, IL 
  (\email{ddoe@imag.com}, \url{http://www.imag.com/\string~ddoe/}).}
\and Paul T. Frank\thanks{Department of Applied Mathematics, Fictional University, Boise, ID 
  (\email{ptfrank@fictional.edu}, \email{jesmith@fictional.edu}).}
\and Jane E. Smith\footnotemark[3]}

\author{
  Boris Martin
  \thanks{Université de Liège, Institut Montefiore (\email{boris.martin@uliege.be}, \email{cgeuzaine@uliege.be})}
  \and
  Pierre Jolivet
  \thanks{Sorbonne Université, CNRS, LIP6 (\email{pierre@joliv.et})} \\
  \and
  Christophe Geuzaine
  \footnotemark[2] \\
}

\usepackage{amsopn}

\ifpdf
\hypersetup{
  pdftitle={Comparison of substructured non-overlapping domain decomposition and
  overlapping additive Schwarz methods for large-scale Helmholtz
  problems with multiple sources },
  pdfauthor={B. Martin, P. Jolivet and C. Geuzaine},
}
\fi

\begin{document}

\maketitle

\begin{abstract}
  Solving large-scale Helmholtz problems discretized with high-order finite elements is notoriously difficult, especially in 3D where direct factorization of the system matrix is very expensive and memory demanding, and robust convergence of iterative methods is difficult to obtain. Domain decomposition methods (DDM) constitute one of the most promising strategy so far, by combining direct and iterative approaches: using direct solvers on overlapping or non-overlapping subdomains, as a preconditioner for a Krylov subspace method on the original Helmholtz system or as an iterative solver on a substructured problem involving field values or Lagrange multipliers on the interfaces between the subdomains. In this work we compare the computational performance of non-overlapping substructured DDM and Optimized Restricted Additive Schwarz (ORAS) preconditioners for solving large-scale Helmholtz problems with multiple sources, as is encountered, e.g., in frequency-domain Full Waveform Inversion. We show on a realistic geophysical test-case that, when appropriately tuned, the non-overlapping methods can reduce the convergence gap sufficiently to significantly outperform the overlapping methods.
\end{abstract}



\section{Introduction}
\label{sec:intro}

Full Waveform Inversion (FWI) is a powerful technique for seismic imaging~\cite{pratt_fwi}, which can estimate subsurface parameters such as velocities by fitting simulated waveform data to observed data. The fitting problem is formulated as an optimization problem, which is typically solved using gradient-based methods, with the gradient computed using adjoint methods. In the frequency-domain, both the misfit evaluation and the adjoint problem require solving a Helmholtz-like PDE with the current subsurface model, e.g., with finite differences or finite elements, resulting in solution of sparse, complex and indefinite linear systems. The size of these linear systems is directly related to desired imaging resolution and thus the frequency: in 3D the size grows at least cubically with frequency. Direct solvers are robust but scale poorly~\cite{mary2017}, although recent progress such as Block Low-Rank compression~\cite{nies_testing_2019, amestoy_obc_blr, amestoy_improving_2013, amestoy_performance_2019, operto_is_2023, mary2017} can somewhat alleviate these issues. When the number of unknowns is on the order of the hundreds of millions, iterative methods become an interesting alternative, despite the challenges of designing efficient preconditioners for indefinite problems~\cite{ErnstGander2012}. Domain decomposition methods (DDM) are a popular choice, leveraging the power of sparse solvers on smaller subdomains. Their main drawback is that, being iterative methods, each RHS requires the full iterative procedure to be applied, making these methods a priori less well-suited for a large number of sources. Nonetheless, DDM has been shown to be able to outperform direct solvers in some cases, even with a fairly large number (130) of sources~\cite{tournier:hal-03942570}.

There exists a wide range of DDMs, with either overlapping or non-overlapping subdomains. One the one hand, the overlap tends to help achieve faster convergence, but at the cost of larger subdomains. On the other hand, non-overlapping methods work on a substructured problem involving a reduced set of unknowns defined on the interfaces between subdomains, resulting in smaller problem sizes. Thus, a non-overlapping approach can theoretically mitigate the computational expense associated with Krylov methods, which require the storage and orthogonalization of a Krylov subspace basis.

Applying an iterative procedure to a large number of sources essentially leads to a cost proportional to that number, but the efficiency of the process can be improved with a parallel implementation that batches large portions of the computations (called \textit{pseudo-block Krylov Methods}), and by using iterative methods that can reduce the total number of iterations for multiple sources. This can be achieved by \textit{recycling} a part of the subspace~\cite{feti2lm} from previous solves, or, when solving for all sources simultaneously, by using a block Krylov method, where information is shared between the different sources~\cite{jolivet_block_iterative_2016}. Essentially, these methods reduce the iteration count, at the expense of additional work per iteration (e.g., additional inner products). Depending on the problem and the context, they may or may not offer a favorable trade-off compared to pseudo-block methods. In particular, the reduction in iteration count is difficult to predict, and may depend on the properties of the matrix and the space spanned by the right-hand sides, the latter being physically related to the choice of source locations~\cite{gabriel2023acceleration}.

In this work, we compare the Optimized Restricted Additive Schwarz (ORAS) method~\cite{oras, oras_variational, bookddm}, which is an overlapping preconditioner for the Helmholtz equation, with the non-overlapping Optimized Schwarz method (OSM)~\cite{Nataf2002, doi:10.1137/S1064827501387012, bookddm, antoine_ddm}, also known as FETI-2LM~\cite{feti2lm}. We focus on analyzing computational times and memory usage for increasingly large 3D problems solved with high-order finite elements on adapted unstructured tetrahedral grids. Meshes are partitioned using METIS~\cite{Karypis_1998} and GMRES serves as Krylov subspace solver. 

The paper is organized as follows. We first introduce the Helmholtz problem and its discretization, as well as our test case of interest in \cref{sec:problem}. The two domain decomposition methods of interest are presented in \cref{sec:ddm}, while the convergence criterion and the parameters for both methods are discussed in ~\cref{sec:parameters}. \Cref{sec:multirhs} then analyzes the results of the numerical experiments in terms of both computational time and memory usage. Conclusions and perspectives for future work are finally given in \cref{sec:conclusion}.

\section{Problem setting}
\label{sec:problem}

\subsection{Helmholtz problem}

We focus on the constant-density, variable-velocity and damping-free Helmholtz equation in a parallelepipedic domain~$\Omega$. For simplicity, we consider a low-order approximation of the Sommerfeld radiation condition, in the form of a Robin boundary condition on the boundary $\Gamma^{\infty}=\partial\Omega$ of $\Omega$. With a point source located at $x=x_s \in \Omega$, the problem writes:
\begin{equation}
  \left\{
    \begin{array}{lll}
    -\Delta u - k^2 u = \delta(x-x_s)  & \text { in } \Omega & \quad (\text{Helmholtz equation}) , \\
    \partial_{\boldsymbol{n}} u - \imath k u=0 & \text { on } \Gamma^{\infty} & \quad (\text{radiation condition}) .
    \end{array}
  \right.
  \label{eq:helmholtz}
\end{equation}
The wavenumber $k$ is a function of space and frequency, and is defined as $k = \frac{\omega}{c(x)}$, where $\omega$ is the angular frequency and $c(x)$ is the velocity field. Our main interest is to solve the equation for multiple source positions $x_s$, which will translate as different right-hand sides in the linear system after finite element discretization. We use a standard (continuous) Galerkin formulation on a tetrahedral mesh adapted to the local wavelength with third-order hierarchical basis functions~\cite{Solin2003}. The resulting linear system is denoted  $A \vec{u} = \vec{f}$, where $\vec{u}$ is the unknown vector, $\vec{f}$ the discretized source, and $A$ is a sparse, complex, symmetric (but not Hermitian) matrix. Except at the lowest frequencies, the matrix is indefinite.

\subsection{A geophysical benchmark: the GO\_3D\_OBS velocity model}
\label{sec:obs}

\begin{figure}
  \centering
  \includegraphics[width=0.8\textwidth]{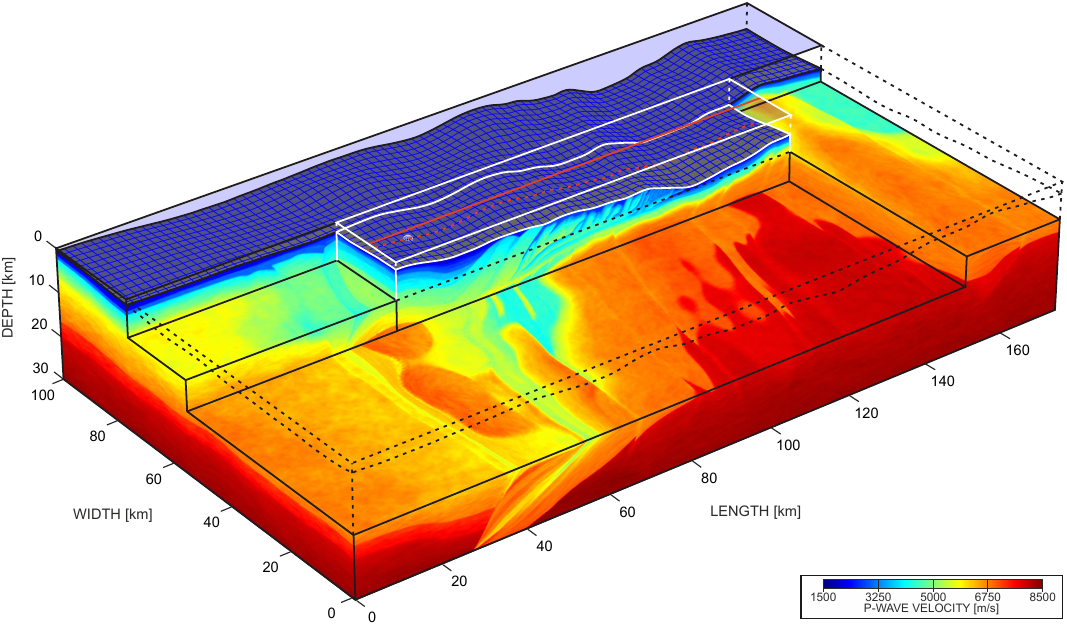}
  \caption{The GO\_3D\_OBS velocity model and the target region~\cite{obs_crustal}.}
  \label{fig:go3d_obs}
\end{figure}

\begin{figure}
  \centering
  \includegraphics[width=0.5\textwidth]{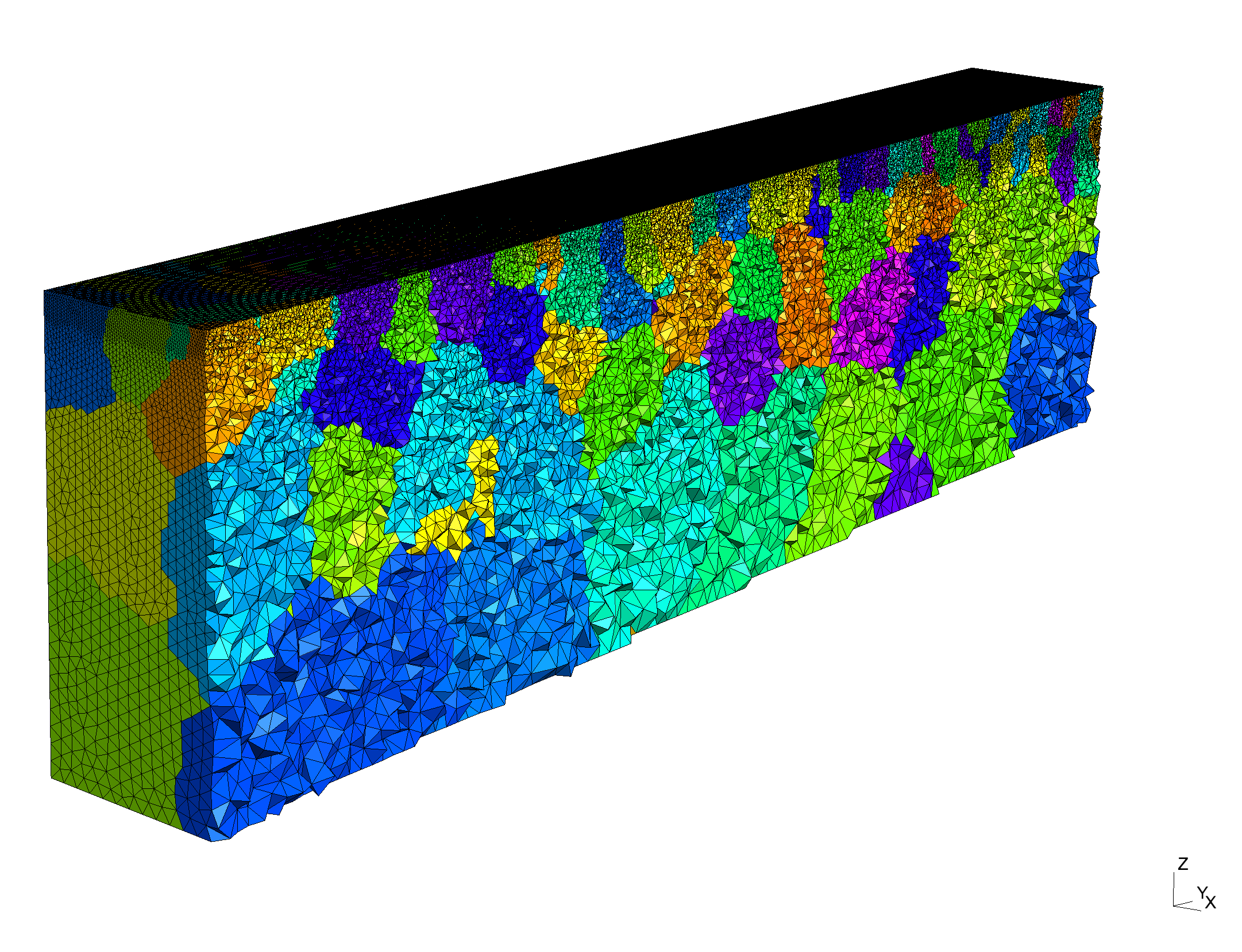}\includegraphics[width=0.5\textwidth]{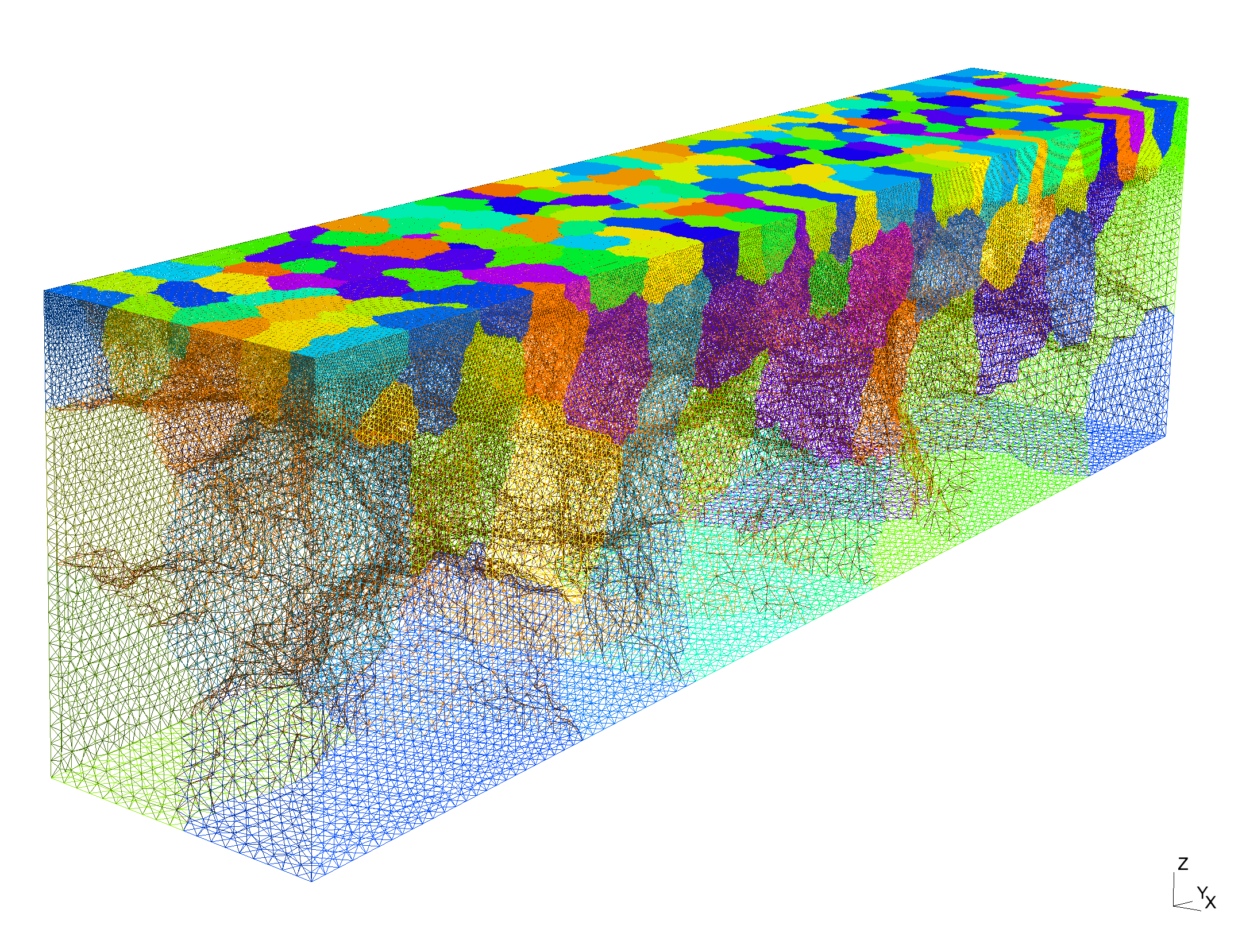}
  \caption{Typical mesh (here with 256 partitions), with element sizes adapted to the local wavelength at \SI{2}{\hertz}.}
  \label{fig:mesh}
\end{figure}

As velocity model, we focus on the GO\_3D\_OBS model~\cite{obs_crustal}, a publicly available 3D geomodel representing a subduction zone, inspired by the geology of the Nankai Trough. The model represents a \qtyproduct{175 x 100 x 30}{\km} domain with a layer of water and complex geological structures. We focus on the \qtyproduct{102 x 20 x 28}{\km} subset of the model, called the target, suggested in \cite{obs_crustal} and illustrated on \cref{fig:go3d_obs}. Velocity ranges from \SI{1500}{\m\per\s} (water) to \SI{8500}{\m\per\s} in some areas. The velocity data is available as point values on a regular grid with a \SI{25}{m} spacing, which is interpolated trilinearly on the nodes of the adapted tetrahedral finite element mesh (cf. \cref{fig:mesh}). From this trilinear interpolation, we project $\frac{1}{c(x)}$ such that $k$ is piecewise linear and continuous on the finite element mesh.

\section{Domain decomposition methods}
\label{sec:ddm}

Let us partition the domain $\Omega$ into $N$ subdomains $\Omega_i$ (for $i=1,..,N$). The subdomains can either be overlapping, or non-overlapping. In both cases, the boundary of $\Omega_i$ is the union of a (possibly empty) portion of $\Gamma^\infty$ and the \textit{artificial} interface $\Sigma_i$ = $\partial \Omega_i \setminus \Gamma^\infty$. If $\mathcal{O}(i)$ denotes the set of indices of the neighbors of subdomain $i$, then $\Sigma_i = \cup_{j \in \mathcal{O}(i)} \Sigma_{ij}$, where:
\begin{itemize}
  \item in the non-overlapping case, $\Sigma_{ij} = \Sigma_{ji}$ is the surface separating $\Omega_i$ and $\Omega_j$;
  \item in the overlapping case, $\Sigma_{ij}$ is the part of $\partial \Omega_i$ inside $\Omega_j$, i.e., $\Sigma_{ij} = \partial \Omega_i \cap \Omega_j$; it is thus distinct from $\Sigma_{ji}$, which lies in $\Omega_i$.
\end{itemize}

On each subdomain, we solve a local problem with impedance-type transmission conditions on the artificial interface:
\begin{equation}
  \left\{
  \begin{array}{lll}
    -\Delta u_i - k^2 u_i = f_i & \text { in } \Omega_{i}, & (\text{Helmholtz equation}) ,\\
    \partial_{\boldsymbol{n_i}} u_i - \imath k u_i = 0 & \text { on } \partial\Omega_i \cap \Gamma^{\infty} & (\text{radiation condition}) ,\\
    \partial_{\boldsymbol{n_i}} u_i - \mathcal{S} u_i = \partial_{\boldsymbol{n_i}} u_j - \mathcal{S} u_j
         & \text { on } \Sigma_{ij}, \ \forall j \in \mathcal{O}(i) & (\text{transmission condition}) ,
  \end{array}
\right.
\label{eq:helmholtz_ddm}
\end{equation}
where $f_i = \delta(x-x_s)$. The operator $\mathcal{S}$ is chosen to improve convergence of DDM. The simplest choice is $\mathcal{S} = \imath k$ \cite{Despres1991}. In this paper, we consider second order operators $\mathcal{S} = \alpha + \beta \Delta_\Sigma$, where the two complex-valued coefficients $\alpha$ and $\beta$ are optimized numerically.

\Cref{eq:helmholtz_ddm} can yield a simple iterative procedure (a Schwarz method) where, at step $m+1$, one can compute $u_i^{(m+1)}$ as a function of $u_j^{(m)}$ for all $j \in \mathcal{O}(i)$. The process is repeated until convergence. In practice, depending on whether the subdomains are overlapping or not, two common algorithms~\cite{bookddm} are the Optimized Restricted Additive Schwarz (ORAS) preconditioner (with overlap), and the substructured solver with unknowns on the interfaces, called hereafter the Optimized Schwarz Method (OSM). We recall both methods in what follows.

\subsection{Overlapping subdomains: the ORAS preconditioner}

From the local solutions $u_i$ on each subdomain, we can compute the global solution $u$ by gluing the local solutions together using a partition of unity. If $E_i(u_i)$ is the extension by zero of $u_i$ on $\Omega$ and $\chi_i(x)$ forms a partition of unity, the global solution can be expressed as:
\begin{equation}
u(x) = \sum_{i=1}^{N} \chi_i(x) E_i(u_i).
\label{eq:partition_unity}
\end{equation}

The ORAS procedure in its continuous form consists of building a sequence $u^{(m)}$ where $u^{(m+1)}$ is computed by restricting $u^{(m)}$ to each subdomain, solving the local problem \cref{eq:helmholtz_ddm} and then summing the results using \cref{eq:partition_unity}. This procedure is a fixed-point iteration that will converge if the transmission conditions are appropriately chosen.

At the discrete level, the restrict, solve, extend, and combine procedure has a natural algebraic equivalent. Let $R_i$ be the restriction matrix that selects the entries related to subdomain $\Omega_i$ from a global vector. $R_i^T$ then extends by zero a local vector to the global space. The discrete partition of unity is a set of diagonal matrices $D_i$ such that $\sum_{i=1}^{N} R_i^TD_iR_i = I$. Typically, the partition is chosen as Boolean, and each degree of freedom (DOF) is \textit{owned} by a single subdomain. Each subdomain solves a problem involving its owned DOFs and some neighboring non-owned DOFs, and outputs new values for the owned ones.

Combining the restrictions, local solves, and extensions weighted by the partition of unity, the ORAS method is expressed as follows:
\begin{equation*}
  \vec{u}^{(k+1)} = \vec{u}^{(k)} + M_{\text{ORAS}}^{-1} (\vec{f}-A\vec{u}^k),
\end{equation*}
with the ORAS preconditioner $M_{\text{ORAS}}^{-1}$ defined by
\begin{equation*}
M_{\text{ORAS}}^{-1} = \sum_{i=1}^{N} R_i^T D_i A_i^{-1} R_i.
\end{equation*}

The matrix $A_i$ is the matrix arising from the discretization of the subproblems \cref{eq:helmholtz_ddm}, with the transmission conditions. This fixed-point iteration can be seen as a preconditioned Richardson iteration, where the preconditioner is the ORAS method.
This preconditioner can be used with other Krylov methods, such as GMRES. In that case, we solve by right-preconditioning:
\begin{equation*}
AM_{\text{ORAS}}^{-1} \vec{y} =  \vec{f},
\end{equation*}
and then $\vec{u} = M_{\text{ORAS}}^{-1}\vec{y}$ is the solution of $A\vec{u}=\vec{f}$.

To summarize, the ORAS method is used as a preconditioner to solve the complete system with GMRES or another Krylov method. The residual involves the FEM matrix $A$ and the preconditioner application involves solving each subproblem once (with $A_i$) and then using the extension and partition of unity to obtain a global vector as output. When using GMRES, the computational cost resides mostly in the preconditioner application, the sparse matrix-vector product and the Gram--Schmidt procedure to build an orthonormal basis of the Krylov subspace.

\subsection{Non-overlapping case: substructured Schwarz method (OSM)}

In the absence of overlap, the interfaces $\Sigma_{ij}$ and $\Sigma_{ji}$ are the same, with opposite normals. On each interface, we define $g_{ji} = \partial_{\boldsymbol{n_i}} u_j - \mathcal{S} u_j$. We note $g$ (without subscript) the tuple of all coupling variables containing, for all subdomains $i$, the $g_{ij}$ corresponding to each subdomain $j$ neighboring $i$.
Using these notations, \cref{eq:helmholtz_ddm} can be rewritten in the form of \cref{eq:helmholtz_ddm_osm}, whose unknowns are the local fields $u_i$ and a tuple of traces $g$:
\begin{equation}
  \left\{
    \begin{array}{lll}
      -\Delta u_i - k^2 u_i = f_i &  \text { in } \Omega_i  & (\text{Helmholtz equation}) \\
      \partial_{\boldsymbol{n_i}} u_i - \imath k u_i = 0 & \text { on } \partial\Omega_i \cap \Gamma^{\infty} & (\text{radiation condition}) \\
      \partial_{\boldsymbol{n_i}} u_i - \mathcal{S} u_i = g_{ji} & \text{ on } \Sigma_{ij}, \ \forall j \in \mathcal{O}(i) & (\text{transmission condition})
    \end{array}
  \right.
\label{eq:helmholtz_ddm_osm}
\end{equation}
At convergence (i.e., for the correct $g$), the solution is continuous, and $u_i = u_j$ on $\Sigma_{ij}$. Combined with the opposite signs of the normals, we obtain the following equation that links the two auxiliary variables (called \textit{interface fields}) and the wavefield, for all pairs of neighboring subdomains:
\begin{equation}
  g_{ij} + g_{ji} = -2\mathcal{S}u \text{ on } \Sigma_{ij}.
  \label{eq:helmholtz_update_g}
\end{equation}

This relation allows us to design a (parallel) Schwarz method that iterates on the interface fields. We first solve, for all $i$, the subdomain \cref{eq:helmholtz_ddm_osm} for the current values of $g_{ji}$ (for all relevant values of $j$), and then, using \cref{eq:helmholtz_update_g}, we compute new values of $g_{ij}$ necessary to solve the problem on subdomain $j$.

\paragraph{Krylov acceleration}
The solution of a local problem \cref{eq:helmholtz_ddm_osm} depends on two sources: $f_i$ and the interface fields $g_{ji}$. By linearity, these can be split as $u_i^{(m)} = v_i + \tilde{u}_i^{(m)}$, where the first part comes from $f_i$ and is constant at each iteration, and the second part represents the wavefield generated by the waves coming from other subdomains. The interface update can be rewritten as the affine fixed-point iteration \cref{eq:helmholtz_affine} for each coupling field:
\begin{equation}
g_{ij}^{(m+1)} = -g_{ji}^{(m)} - 2\mathcal{S} \tilde{u}_i^{(m)} - 2\mathcal{S} v_i.
\label{eq:helmholtz_affine}
\end{equation}

Combining this for each interface equation, the global update can be written as $g^{(m+1)} = \mathcal{T}g^{(m)} + b$ where $b$ represents the outgoing waves generated by all the sources $f_i$ and $\mathcal{T}$ is an exchange operator acting on the tuple of coupling variables~\cite{bookddm}. This can be reformulated as the solution of the following substructured problem in terms of $g$:
\begin{equation}
  (\mathcal{I}-\mathcal{T})g = b.
  \label{eq:osm_substructured}
\end{equation}
The vector $b$ is obtained by solving the local problems with homogeneous boundary conditions. Then, to apply $\mathcal{T}$, one solves the subproblems without $f_i$ and with inhomogeneous boundary conditions.

At the discrete level, this system can be solved with a Krylov method such as GMRES (and the fixed-point iteration is actually the Richardson iteration applied to it). Once the interface fields are recovered, the complete wavefield can be computed by solving the local problems with both these boundary conditions and the local source. It is worth noting that \cref{eq:osm_substructured} can also be seen as a Schur complement: combining \cref{eq:helmholtz_ddm_osm} for all subdomains and \cref{eq:helmholtz_update_g} for all interfaces yields a global system, and the substructured problem is actually its Schur complement with respect to all the $u_i$, and the fixed-point iteration is a block-Jacobi iteration applied to the Schur complement.

In this approach, GMRES is applied to the interface problem, without preconditioner, but with an operator that contains local solves. In this workflow, the cost is spread between the local solves, the assembly of inhomogeneous terms and the construction of an orthonormal basis of reduced dimension of the Krylov subspace.

\subsection{A priori comparison of ORAS and OSM}

Both ORAS and OSM are based on the principle of solving subproblems with transmission boundary conditions until the full solution is recovered. Thus, the main computational cost in both cases will be to factorize the local matrices (once) and to solve the resulting triangular systems at each iteration. The workflow is however different: ORAS is a preconditioner for solving the full problem iteratively, whereas in OSM the local solves are part of the operator itself, and no preconditioner is used.
The OSM method should a priori be more economical for several reasons:
\begin{itemize}
\item The absence of overlap reduces the subdomain size, making the local solves cheaper, and the factorization less memory-intensive.
\item The iterative method is applied on the interface fields, which leads to fewer unknowns (typically between 5 and 10 times less) than the full problem. This reduces the cost (in both time and memory) of the Gram--Schmidt orthogonalization procedure in GMRES.
\item The global sparse matrix-vector product disappears. However, it is replaced by a sparse matrix-vector product (SpMV) on the local interface fields to compute the right-hand side of the local problems at the next iteration. Again, this should be an order of magnitude cheaper.
\end{itemize}

However, ORAS can converge faster thanks to the overlap, and is more flexible as it is a preconditioner. In particular, inexact subdomain solves can be used safely: they can slow down convergence but not yield an incorrect result. Conversely, an inexact solve in OSM will impact the resulting solution as it introduces an error in the operator instead of the preconditioner.

In the following sections, we will evaluate how the convergence rates differ between the two methods. Afterwards, we will run experiments on large problems with multiple sources to assess the efficiency of the two methods in terms of time and memory consumption in realistic application settings.

\subsection{Implementation}

For our experiments, we used Gmsh~\cite{Geuzaine2009Gmsh} to generate tetrahedral meshes, which are then partitioned using METIS~\cite{Karypis_1998}. We extended Gmsh to generate the overlaps needed for ORAS. The problem assembly is done with GmshFEM and GmshDDM~\cite{ORBi-fem, ORBi-royer-these}, a Gmsh-based FEM solver and its extension designed for non-overlapping domain decomposition methods. For both ORAS and OSM, GMRES is used through PETSc~\cite{petsc-user-ref} and its binding to HPDDM~\cite{Jolivet2021KSPHPDDM}, which provides a pseudo-block GMRES implementation suited for multiple RHSs. Both approaches use MUMPS~\cite{mumps1,mumps2} as a solver (called through PETSc) to compute a $LDL^T$ factorization of the subproblem matrices and solve the resulting triangular problems.

While ORAS with multiple RHS is straightforward with these tools (as the additive Schwarz preconditioner of PETSc is already adapted for multiple RHS), OSM is more complex. GmshDDM was first implemented as a custom matrix-vector product provided to PETSc. We upgraded it to handle matrix-matrix products, i.e., applying the $\mathcal{I}-\mathcal{T}$ operator to a batch of vectors. This required reformulating the artificial sources assembly as a mass-like linear operator that can be precomputed then applied to a batch of incoming waves efficiently. Finally, in OSM, the local problems have a block-triangular $2 \times 2$ structure, as the local wavefield can be computed first, and the updated outgoing wavefield depends on it. Solving directly this problem is inefficient, as we lose the symmetry and have a larger system. Exploiting the triangular structure allows us to take advantage of the symmetry of the two diagonal blocks (the Helmholtz operator and a mass-matrix on the interface), and can be done through PETSc with a \textit{fieldsplit} preconditioner. We extended this preconditioner to handle efficiently the multiple RHS case, as the current version of PETSc simply loops over all of these. More details are available in \cref{sec:reproducibility}.

\section{Choice of parameters and convergence criterion}
\label{sec:parameters}

It is well-known that appropriate transmission conditions are critical to achieve good convergence in DDM, both without overlap~\cite{dolean2017optimized, doi:10.1137/S1064827501387012} and with overlap~\cite{doi:10.1137/15M1021659}. The simplest choice is $\mathcal S =\imath k$ but more complex conditions can be chosen: second-order operators using the surface Laplacian, Perfectly Matched Layers~\cite{bootland2023numericalassessmentpmltransmission, royer2022_pml} or high-order rational functions~\cite{bounbendir2012_pade, modave2020_pade}. We choose to use second-order operators of the form $\mathcal S = \alpha + \beta \Delta_\Sigma$, which provide a good compromise between cost and performance. The complex coefficients $\alpha$ and $\beta$ should be chosen to balance the convergence rate of all the wave modes (propagating, grazing, evanescent). In this work we numerically optimized them to find the coefficients leading to the minimal number of Krylov iterations. In addition, for ORAS we always choose the minimal overlap size, i.e., each subdomain for ORAS is constructed by adding one layer of elements (connected by at least one node) to the elements of the corresponding OSM subdomain.

Since ORAS and OSM have different residuals (preconditioned or unpreconditioned residual from the matrix in ORAS, residual on the interface problem in OSM), we need to define a common convergence criterion. The evolution of the $L^2$ error over iterations is measured as well as the residuals, to find stopping criteria that yield similar accuracies with both approaches. With ORAS, the unpreconditioned residual (from the right-preconditioned GMRES) is used, as all experiments showed that right-preconditioned GMRES was slightly faster, and required a looser tolerance. In this section, we perform our experiments with double precision scalars. Single precision scalars give similar results until we lower the stopping residual below  about $10^{-4}$.

\subsection{Homogeneous medium}

We begin by simulating the \qtyproduct{102 x 20 x 28}{\km} target from GO\_3D\_OBS model by considering that it is entirely made of water ($c = \SI{1500}{\m\per\s}$), with a source at \SI{1}{\hertz} positioned near the surface and impedance conditions applied on the 6 faces of the box.

The box is discretized with an unstructured tetrahedral mesh with prescribed uniform mesh density of 4 points per wavelength. With third-order finite element basis functions this leads to around 20 million unknowns. The METIS graph partitioner is used to create 256 subdomains.

\begin{figure}
    \begin{tikzpicture}
    \begin{axis}[
        width=12cm,height=6cm,
        xlabel={Iteration},
        ylabel={Relative $L^{2}$ error},
        ymode=log,
        legend style={draw=none,fill=none,
        at={(0.02,0.02)}, anchor=south west},
        legend cell align=left,
    ]
    \addplot+[mark=o, mark repeat=5] table[
            col sep=semicolon,
            y index=0,
            x expr=\coordindex+1 ] {conv_final_homog_metis_double_l2_error_ddm_field_3_g_3_oo0_alpha_0.0_pi_ndom_256.csv};
    \addlegendentry{OSM -- Order 0}

    \addplot+[mark=triangle, mark repeat=5] table[
            col sep=semicolon,
            y index=0,
            x expr=\coordindex+1
        ] {conv_final_homog_metis_double_l2_error_ddm_field_3_g_3_oo2_alpha_0.5_pi_ndom_256.csv};
    \addlegendentry{OSM -- Order 2}

    \addplot+[mark=square, mark repeat=5] table[
            col sep=semicolon,
            y index=0,
            x expr=\coordindex+1
        ] {conv_final_homog_metis_double_l2_error_oras_oo0_alpha_0.0_pi_ndom_256.csv};
    \addlegendentry{ORAS -- Order 0}

    \addplot+[mark=diamond, mark repeat=5] table[
            col sep=semicolon,
            y index=0,
            x expr=\coordindex+1
        ] {conv_final_homog_metis_double_l2_error_oras_oo2_alpha_0.5_pi_ndom_256.csv};
    \addlegendentry{ORAS -- Order 2}
    \end{axis}
    \end{tikzpicture}
    \caption{Convergence of the $L^2$ error for the homogeneous medium with 256 subdomains.}
    \label{fig:conv_256_homoh}
    \end{figure}
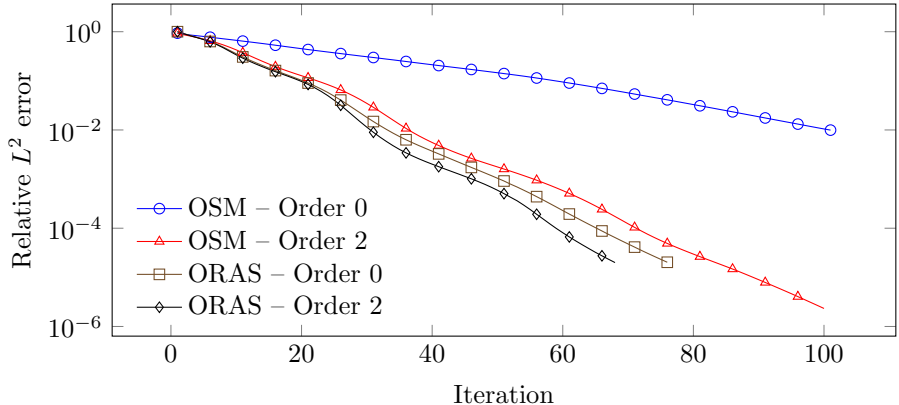

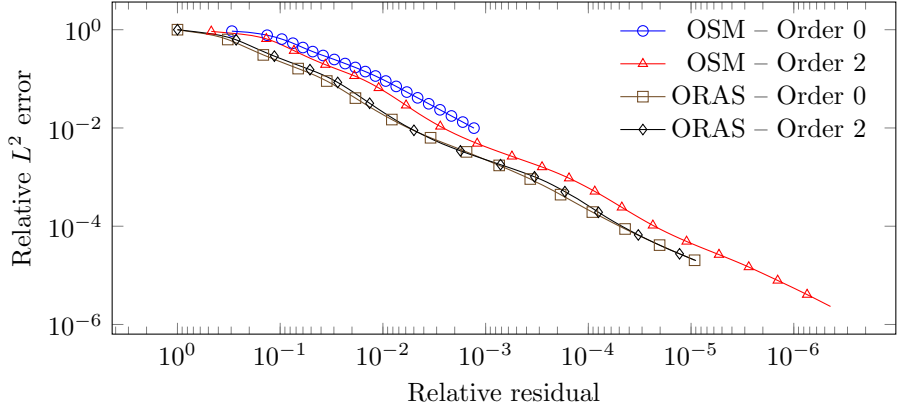
\begin{figure}
\begin{tikzpicture}
        \begin{axis}[
            width=12cm,height=6cm,
            xlabel={Relative residual},
            ylabel={Relative $L^{2}$ error},
            xmode=log, ymode=log,
            x dir=reverse,
            legend style={
                at={(0.98,0.98)},    anchor=north east,   draw=none,fill=none
            },
            legend cell align=right,
        ]
        \addplot+[mark=o, mark repeat=5] table[
                col sep=semicolon,
                x expr=\thisrowno{1} / 2.28738e-05, y index=0  ] {conv_final_homog_metis_double_l2_error_ddm_field_3_g_3_oo0_alpha_0.0_pi_ndom_256.csv};
        \addlegendentry{OSM -- Order 0}

        \addplot+[mark=triangle, mark repeat=5] table[
                col sep=semicolon,
                x expr=\thisrowno{1} / 9.58934e-05, y index=0
            ] {conv_final_homog_metis_double_l2_error_ddm_field_3_g_3_oo2_alpha_0.5_pi_ndom_256.csv};
        \addlegendentry{OSM -- Order 2}

        \addplot+[mark=square, mark repeat=5] table[
                col sep=semicolon,
                x index=1, y index=0
            ] {conv_final_homog_metis_double_l2_error_oras_oo0_alpha_0.0_pi_ndom_256.csv};
        \addlegendentry{ORAS -- Order 0}

        \addplot+[mark=diamond, mark repeat=5] table[
                col sep=semicolon,
                x index=1, y index=0
            ] {conv_final_homog_metis_double_l2_error_oras_oo2_alpha_0.5_pi_ndom_256.csv};
        \addlegendentry{ORAS -- Order 2}
        \end{axis}
        \end{tikzpicture}
        \caption{Relationship between the relative residual and the relative $L^2$ error for the homogeneous medium with 256 subdomains.}
        \label{fig:conv_256_residual_homog}
    \end{figure}

\Cref{fig:conv_256_homoh} shows the convergence of the $L^2$ error for both methods with different transmission conditions: the simplest zeroth-order condition with $\beta=0$ (``Order 0'') and the optimized second order condition (``Order 2''). Clearly, ORAS has a better convergence than OSM, whether the transmission conditions are optimized or not. However, optimizing them dramatically reduces the gap, making OSM almost as effective. \Cref{fig:conv_256_residual_homog} shows the relationship between the relative error and the relative residual. We observe that, for a target error, the OSM tolerance should be stricter than with ORAS.

\subsection{Heterogeneous medium}

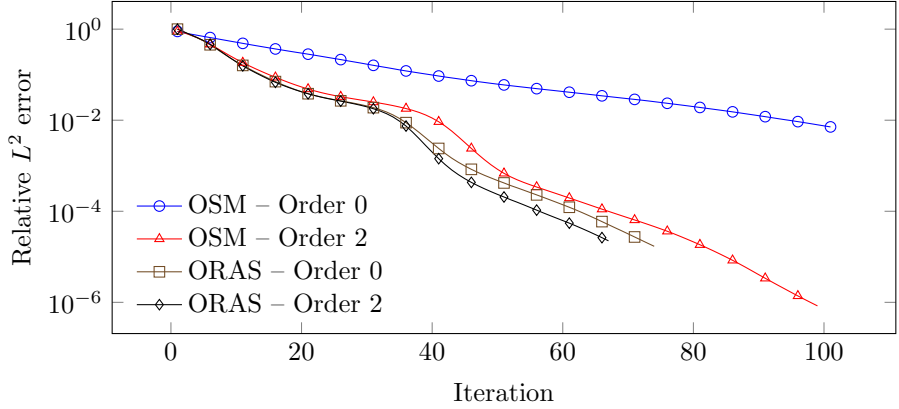
\begin{figure}
\begin{tikzpicture}
\begin{axis}[
    width=12cm,height=6cm,
    xlabel={Iteration},
    ylabel={Relative $L^{2}$ error},
    ymode=log,
    legend style={draw=none,fill=none,
        at={(0.02,0.02)}, anchor=south west},
    legend cell align=left,
]
\addplot+[mark=o, mark repeat=5] table[
        col sep=semicolon,
        y index=0,
        x expr=\coordindex+1 ] {conv_final_double_l2_error_ddm_field_3_g_3_oo0_alpha_0.0_pi_ndom_256.csv};
\addlegendentry{OSM -- Order 0}

\addplot+[mark=triangle, mark repeat=5] table[
        col sep=semicolon,
        y index=0,
        x expr=\coordindex+1
    ] {conv_final_double_l2_error_ddm_field_3_g_3_oo2_alpha_0.5_pi_ndom_256.csv};
\addlegendentry{OSM -- Order 2}

\addplot+[mark=square, mark repeat=5] table[
        col sep=semicolon,
        y index=0,
        x expr=\coordindex+1
    ] {conv_final_double_l2_error_oras_oo0_alpha_0.0_pi_ndom_256.csv};
\addlegendentry{ORAS -- Order 0}

\addplot+[mark=diamond, mark repeat=5] table[
        col sep=semicolon,
        y index=0,
        x expr=\coordindex+1
    ] {conv_final_double_l2_error_oras_oo2_alpha_0.5_pi_ndom_256.csv};
\addlegendentry{ORAS -- Order 2}
\end{axis}
\end{tikzpicture}
\caption{Convergence of the $L^2$ error for the homogeneous medium with 256 subdomains.}
\label{fig:conv_256}
\end{figure}

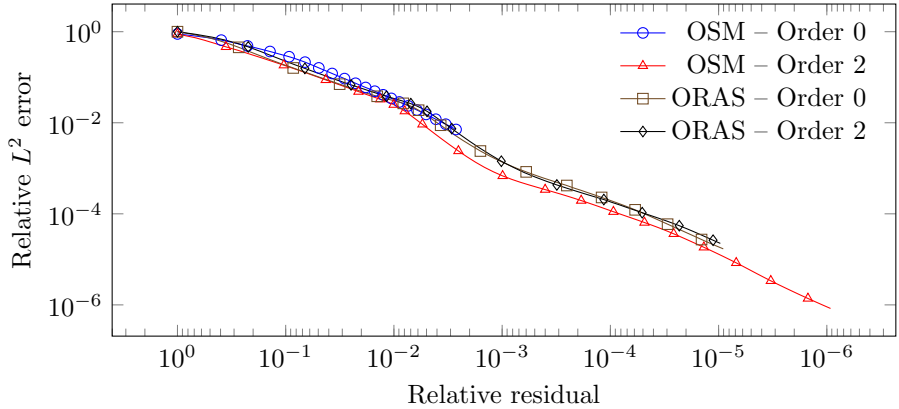
\begin{figure}
\begin{tikzpicture}
    \begin{axis}[
        width=12cm,height=6cm,
        xlabel={Relative residual},
        ylabel={Relative $L^{2}$ error},
        xmode=log, ymode=log,
        x dir=reverse,
        legend style={
            at={(0.98,0.98)},    anchor=north east,   draw=none,fill=none
        },
        legend cell align=right,
    ]
    \addplot+[mark=o, mark repeat=5] table[
            col sep=semicolon,
            x expr=\thisrowno{1} / 2.28738e-05, y index=0  ] {conv_final_double_l2_error_ddm_field_3_g_3_oo0_alpha_0.0_pi_ndom_256.csv};
    \addlegendentry{OSM -- Order 0}

    \addplot+[mark=triangle, mark repeat=5] table[
            col sep=semicolon,
            x expr=\thisrowno{1} / 9.58934e-05, y index=0
        ] {conv_final_double_l2_error_ddm_field_3_g_3_oo2_alpha_0.5_pi_ndom_256.csv};
    \addlegendentry{OSM -- Order 2}

    \addplot+[mark=square, mark repeat=5] table[
            col sep=semicolon,
            x index=1, y index=0
        ] {conv_final_double_l2_error_oras_oo0_alpha_0.0_pi_ndom_256.csv};
    \addlegendentry{ORAS -- Order 0}

    \addplot+[mark=diamond, mark repeat=5] table[
            col sep=semicolon,
            x index=1, y index=0
        ] {conv_final_double_l2_error_oras_oo2_alpha_0.5_pi_ndom_256.csv};
    \addlegendentry{ORAS -- Order 2}
    \end{axis}
    \end{tikzpicture}
    \caption{Relationship between the relative residual and the relative $L^2$ error for the homogeneous medium with 256 subdomains.}
    \label{fig:conv_256_residual}
\end{figure}

We now consider the same geometry but with the velocity field from the GO\_3D\_OBS dataset. The mesh size is locally adapted to keep 4 points per wavelength. The adaptivity allows us to use fewer elements than in the homogeneous case, as the wavelength is around 5 times smaller in water than in soil. In order to keep a problem of similar size, the frequency is increased to \SI{2}{\hertz}, which again leads to about 20 million unknowns. We again look at the $L^2$ convergence (\cref{fig:conv_256}) as well as the relationship between the $L^2$ error and the residuals (\cref{fig:conv_256_residual}). While ORAS exhibits again superior convergence, using optimized conditions makes the difference smaller. Compared to the homogeneous case, the residual/error curves are much closer, and similar tolerances yield similar accuracies. In what follows, we thus use the same tolerance for both methods.

\section{Comparison on large-scale problems with multiple right-hand sides}
\label{sec:multirhs}

When solving multiple sources, the setup cost (matrix assembly and factorization) is shared, whereas the iteration cost is duplicated. On paper, the cost for solving $N_s$ sources should be proportional to $N_s$ for the iterative part, but there is a large potential for batching. In particular, the subdomain triangular solves can handle multiple right-hand sides efficiently using vectorized operations (replacing level 2 BLAS calls with level 3 BLAS) and better cache locality, as the factors must be loaded once instead of $N_s$ times. As shown in~\cite{jolivet_block_iterative_2016}, having larger groups of RHS can be more efficient, especially when the solve is multi-threaded.

For these reasons, even when block Krylov methods or subspace recycling are not used, a good treatment of the multiple RHS is critical. When $N_s$ is large, one can either solve the $N_s$ RHS simultaneously, solve each of them sequentially, or split them into smaller groups (e.g., 8 or 16 at once), called the \textit{batch size}. The larger the batch size, the more efficient the solve, at the expense of a higher memory consumption, as more vectors must be stored in the Krylov method. Careful trade-offs are necessary, as using a large batch size might consume memory that could be better allocated to other purposes, such as increasing subdomain sizes or setting a higher GMRES restart parameter.

Overall, both memory and computational time should be considered. In what follows, we focus on minimizing the computational time for a given hardware, which limits the available memory. We first discuss the memory consumption of both methods to evaluate how large the problems we solve can be, given the available hardware and implementation. Then, we compare the computational times of the two methods for increasingly larger problems, with different batch sizes.

All experiments were run on the CPU partition of the LUMI supercomputer\footnote{\url{https://www.lumi-supercomputer.eu/}}. Each node consists of 2 AMD EPYC 7763 CPUs, each with 64 cores, so that a node has 128 cores and \SI{256}{\giga\byte} of RAM. In each case, we used a tolerance of $10^{-4}$ on the relative residual, and a GMRES restart parameter of 50. In this case, we use single precision scalars, which are sufficient to reach the stopping criterion and allow large savings.

\subsection{Memory analysis of a typical case}

Besides fixed costs such as mesh storage and indexing, the memory cost of a simulation is primarily dominated by the factorization and the subspace storage in GMRES. When using a large number of small subdomains, a reduced amount of memory is required compared to using few large subdomains. However, this reduction in memory usage comes with a slower convergence rate. Memory usage from matrix factorizations depends essentially on the subdomain size, whereas the cost of GMRES can be controlled through the batch size and the restart parameter. The latter can however impact convergence, and smaller batches may limit the algorithm efficiency. 

The Krylov subspace to store is significantly smaller with OSM than with ORAS, as it contains vectors representing interface data rather than volume data. This efficiency can be leveraged by using the freed resources to increase the batch size, the restart parameter, or the problem size. It is however worth noting that this interface problem grows larger as the number of partitions increases, and the difference with ORAS can become negligible for extremely small subdomain sizes.

\paragraph{Numerical experiments}

To estimate orders of magnitude, we measure memory consumption for a \SI{2}{\hertz} simulation of our test model, with 64 sources. This leads to around 20 million complex-valued unknowns. Since the LUMI cluster has nodes with \SI{256}{\giga\byte} and 128 cores, we aim for a partitioning where each subdomain can be solved on a single-thread process and less than \SI{2}{\giga\byte}. We partition into either 256, 512, or 1024 subdomains and run the simulations accordingly.

\begin{table}[ht]
  \centering
  \begin{tabular}{|c|c|c|c|c|c|c|c|}
       \hline
\multirow{2}{*}{{Method}}
       & \multirow{2}{*}{$N$}
       & \multicolumn{3}{c|}{{Factorization (MiB)}}
       & \multicolumn{3}{c|}{{Factors only (MiB)}} \\
       \cline{3-8}
& & {Mean} & {Min} & {Max}
         & {Mean} & {Min} & {Max} \\
       \hline
\multirow{3}{*}{ORAS} & 256 & 732 & 544 & 899 & 497 & 381 & 594 \\
\cline{2-8}
       & 512 & 341 & 262 & 417 & 226 & 171 & 277 \\
\cline{2-8}
       & 1024 & 163 & 108 & 206 & 106 & 70 & 128 \\
\hline
\multirow{3}{*}{OSM} & 256 & 439 & 346 & 494 & 305 & 267 & 338 \\
\cline{2-8}
       & 512 & 175 & 139 & 200 & 120 & 98 & 133 \\
\cline{2-8}
       & 1024 & 71 & 52 & 86 & 48 & 38 & 54 \\
  \hline
  \end{tabular}
  \caption{Memory usage for factorization, and for storing the factors (i.e., after freeing the intermediate buffers) for a 20M unknowns problem. Values are provided as minimum, maximum, and average over all processes, with each process managing a subdomain.}
  \label{tab:memory_peak_factor}
\end{table}

\begin{table}[ht]
  \centering
  \begin{tabular}{|c|c|c|c|c|c|}
       \hline
\multirow{2}{*}{{Method}}
       & \multirow{2}{*}{$N$}
       & \multicolumn{2}{c|}{{Local vector size}}
       & \multicolumn{2}{c|}{{Memory for 64 RHSs (MiB)}} \\
       \cline{3-6}
& & {Min} & {Max}
         & {Min} & {Max}  \\
       \hline
\multirow{3}{*}{ORAS} & 256 & 74,701 & 88,621 & 1824 & 2164 \\
\cline{2-6}
       & 512 & 36,490 & 45,013 & 891 & 1099 \\
\cline{2-6}
       & 1024 & 17,736 & 23,296 & 433 & 569 \\
\hline
\multirow{3}{*}{OSM} & 256 & 4,923 & 15,129 & 120 & 369 \\
\cline{2-6}
       & 512 & 3,226 & 10,413 & 79 & 254 \\
\cline{2-6}
       & 1024 & 2,019 & 6,832 & 49 & 167 \\
  \hline
  \end{tabular}
  \caption{Local size of the vector over which GMRES iterates, and memory needed to store a full batch (i.e., $64 \times 50$ vectors with single-precision complex scalars) for the 20M unknowns problem.}
  \label{tab:memory_system_batch}
\end{table}

\begin{table}[ht]
  \centering
  \begin{tabular}{|c|c|c|c|c|}
       \hline
\multirow{2}{*}{{Method}}
       & \multirow{2}{*}{$N$}
       & \multicolumn{3}{c|}{{Peak memory per process (MiB)}} \\
       \cline{3-5}
& & {Mean} & {Min} & {Max}  \\
  \hline
\multirow{3}{*}{ORAS} & 256 & 3703 & 3489 & 4003 \\
\cline{2-5}
       & 512 & 2014 & 1839 & 2398 \\
\cline{2-5}
       & 1024 & 1336 & 1024 & 1706 \\
\hline
\multirow{3}{*}{OSM} & 256 & 1532 & 1324 & 1727 \\
\cline{2-5}
       & 512 & 1032 & 1026 & 1089 \\
\cline{2-5}
       & 1024 & 1027 & 1022 & 1045 \\
  \hline
  \end{tabular}
  \caption{Maximum resident set size (RSS) per process, measured at runtime on our implementation.}
  \label{tab:memory_total}
\end{table}

\Cref{tab:memory_peak_factor} reports for each partitioning and either OSM or ORAS the memory needed for matrix factorization, both as peak consumption and to store the resulting coefficients. In the absence of overlap, the factorization is significantly cheaper: depending on the subdomain sizes, the number of overlapping degrees of freedom ranges from 30\% to 50\% of the total number of DOFs in a subdomain, despite the overlap being chosen as small as possible. In both OSM and ORAS, a finer partitioning leads to smaller factorization costs, in a way which is consistent with the known results that, for a 3D problem of size $N$, memory scales asymptotically as $N^{\frac{4}{3}}$~\cite{mary2017}: every time we double the number of subdomains, the memory cost is reduced by a factor larger than 2. Overall, we see that large subdomains can be quite expensive, especially with ORAS: the partition in 256 subdomains leads us to having a peak memory of \SI{899}{\mebi\byte}, which is almost half the budget of one process.

\Cref{tab:memory_system_batch} reports the local size of the vector over which GMRES iterates (local part of the full vector in ORAS, local interface fields in OSM) as well as the memory needed to do a full batch simulation of 64 sources and a restart parameter of 50. Scalars are stored as single-precision complex number, taking 8 bytes. It highlights another perk of using OSM: since GMRES iterates on the interface fields, the Krylov subspace has a reduced size. This again leads to a smaller memory footprint, especially for larger subdomains. As smaller subdomains lead to more interfaces, the \textit{total} number of interface DOFs increases in OSM when the partitioning gets finer. Thus, the memory advantage of OSM is especially attractive for large subdomains, as the memory saving can be a factor larger than 6. It is however worth noting that interface DOFs suffer from a suboptimal load balancing, as the number of interface unknowns can vary significantly between subdomains due to their irregular shapes. As the table shows, the ratio between the subdomain with the most and the least interface unknowns can reach 3. One could optimize memory usage by distributing domains among nodes such that each node has some domains with many interfaces and some with fewer interfaces.

\pgfplotsset{compat=1.18}
\usepgfplotslibrary{groupplots}
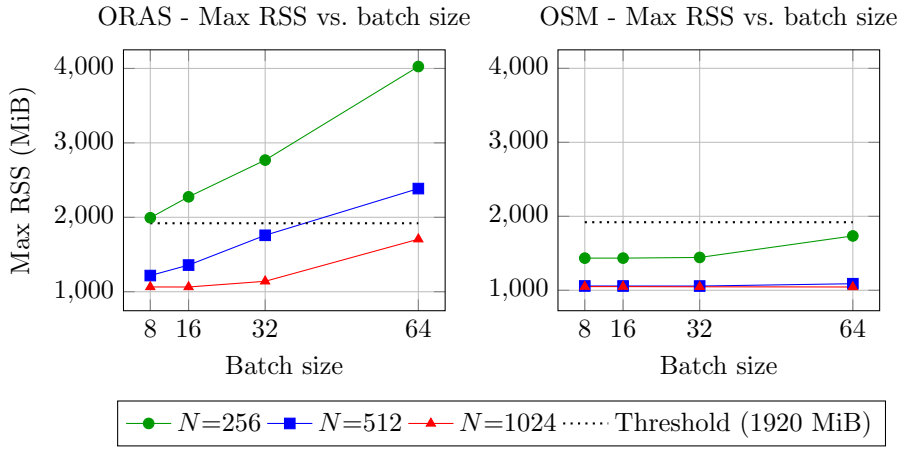
\begin{figure}[ht]
  \centering
  \begin{tikzpicture}
    \begin{groupplot}[
      group style={group size=2 by 1, horizontal sep=1.5cm, ylabels at=edge left},
      width=0.45\textwidth,
      xlabel={Batch size},
      ylabel={Max RSS (MiB)},
      grid=both,
      xtick={8,16,32,64},
      legend to name=sharedlegend,
      legend style={at={(0.5,-0.15)}, anchor=north, legend columns=-1},
      every mark/.append style={solid},
      legend image post style={mark options={solid}},
      ymax=4250
    ]

    \nextgroupplot[title={ORAS - Max RSS vs.\ batch size}]
    \addplot+[mark=*, color=green!60!black] coordinates { (8,1993.878906) (16,2275.160156) (32,2768.664062) (64,4025.417969) };
    \addlegendentry{$N$=256}
    \addplot+[mark=square*, color=blue] coordinates { (8,1217.328125) (16,1358.097656) (32,1758.230469) (64,2385.750000) };
    \addlegendentry{$N$=512}
    \addplot+[mark=triangle*, color=red] coordinates { (8,1064.066406) (16,1063.539062) (32,1140.218750) (64,1705.636719) };
    \addlegendentry{$N$=1024}
    \addplot[dotted, thick, color=black] coordinates {(8, 1920) (64, 1920)};
    \addlegendentry{Threshold (1920 MiB)}
    \nextgroupplot[title={OSM - Max RSS vs.\ batch size}]
    \addplot+[mark=*, color=green!60!black] coordinates { (8,1434.578125) (16,1434.828125) (32,1443.718750) (64,1734.695312) };
    \addlegendentry{$N$=256}
    \addplot+[mark=square*, color=blue] coordinates { (8,1058.199219) (16,1057.679688) (32,1057.156250) (64,1088.585938) };
    \addlegendentry{$N$=512}
    \addplot+[mark=triangle*, color=red] coordinates { (8,1049.750000) (16,1049.714844) (32,1047.343750) (64,1044.429688) };
    \addlegendentry{$N$=1024}
    \addplot[dotted, thick, color=black] coordinates {(8, 1920) (64, 1920)};
    \addlegendentry{Threshold (1920 MiB)}
    \end{groupplot}
    \node at ($(group c1r1.south)!0.5!(group c2r1.south) - (0,1.5cm)$) {\pgfplotslegendfromname{sharedlegend}};
  \end{tikzpicture}
  \caption{Maximum resident set size (RSS) over all processes, vs.\ batch size for OSM and ORAS, grouped by number of subdomains $N$.}
\label{fig:memory_2hz_batch}

\end{figure}
 
Finally, \cref{tab:memory_total} presents the total memory usage at runtime. It includes not only the factors and GMRES space discussed above, but also all quantities we neglected beforehand: intermediary buffers, mesh data, the matrices, the velocity model, etc. With OSM, large subdomains easily fit into one process (\SI{2}{\giga\byte}), unlike in the ORAS case. This issue can be mitigated by reducing the batch size or the restart parameter. However, even without a large Krylov subspace to store, the ORAS case with only 256 subdomains remains too large to fit in single-thread MPI processes. The 512 subdomains case almost fits on average, and can be made to fit using a batch size of 32 instead of 64. We will thus focus on cases with either two threads per process, or a fine enough partitioning when using ORAS. \Cref{fig:memory_2hz_batch} showcases different values of the maximum memory (over all processes) for different batch sizes, with a red-dotted line representing the available memory per process.

Overall, we see that OSM is substantially more efficient in terms of memory, thanks to both the substructuring and the smaller subdomain size.
It allows us to use larger batch sizes at a reasonable cost and to solve larger problems for a given hardware. In particular, we can solve problems with 20 million unknowns using only 2 nodes of the LUMI cluster, leading to one node per 10 million unknowns. If we keep the subdomain size constant, this means huge problems (over a billion unknowns) could be solved using only a fraction of the cluster.

\subsection{Computational time analysis on large cases}

\paragraph{Influence of the batch size}
We first compare different batch sizes for a \SI{4}{\hertz} problem (165 million unknowns) on 4096 subdomains and 2 threads per MPI process, which ensures enough memory is available. As \cref{fig:batch_size_time_per_it} shows, increasing the batch size leads to faster computations on all aspects. However, moving from 32 to 64 does not yield additional speedup. This is consistent with findings from~\cite{tournier:hal-03942570} where the authors mention an optimal size around 20.

\pgfplotsset{compat=1.17}
\begin{figure}
  \centering
\begin{tikzpicture}
  \begin{groupplot}[
    group style={
      group size=2 by 1,
      horizontal sep=1.5em,
      ylabels at=edge left,
      xlabels at=edge bottom,
    },
    scale only axis,
width=0.4\textwidth,
ymin=0,
    ymax=4,
    xlabel={Batch size},
    ylabel={Time per iteration for 64 sources (s)},
    xmode=log,
    log basis x=2,
    xtick={2,4,8,16,32,64},
    grid=major,
  ]

\nextgroupplot[title={ORAS},
    legend to name=sharedlegend,
    legend columns=3,]
      \addplot[blue,        mark=*,         thick]
        table[col sep=comma,x=batch,y expr=\thisrow{tgmres_it}]
          {multirhs_post_batch_4096_out_ORAS.csv};
      \addlegendentry{GMRES orthogonalization}
      \addplot[red,         mark=square*,   thick]
        table[col sep=comma,x=batch,y expr=\thisrow{tmvp_it}]
          {multirhs_post_batch_4096_out_ORAS.csv};
      \addlegendentry{SpMM}
      \addplot[green!60!black,mark=triangle*,thick]
        table[col sep=comma,x=batch,y expr=\thisrow{tmatsolve_it}]
          {multirhs_post_batch_4096_out_ORAS.csv};
      \addlegendentry{Local solves}

\nextgroupplot[
      title={OSM},
    ]
      \addplot[blue,        mark=*,         thick]
        table[col sep=comma,x=batch,y expr=\thisrow{tgmres_it}]
          {multirhs_post_batch_4096_out_OSM.csv};
      \addplot[red,         mark=square*,   thick]
        table[col sep=comma,x=batch,y expr=\thisrow{tmvp_it}]
          {multirhs_post_batch_4096_out_OSM.csv};
      \addplot[green!60!black,mark=triangle*,thick]
        table[col sep=comma,x=batch,y expr=\thisrow{tmatsolve_it}]
          {multirhs_post_batch_4096_out_OSM.csv};

  \end{groupplot}

\end{tikzpicture}
\bigskip
\centering
\pgfplotslegendfromname{sharedlegend}
\caption{Time per iteration (normalized for full batch) for ORAS and OSM, for the \SI{4}{\hertz} problem and 4096 subdomains.}
\label{fig:batch_size_time_per_it}
\end{figure}
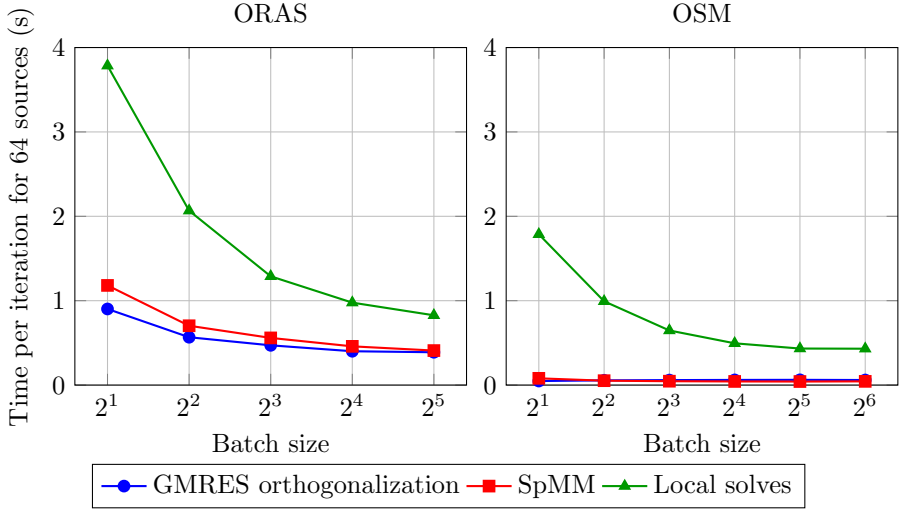
 
\paragraph{Weak scaling}
We evaluate the scaling of both methods on larger cases, using similar domain sizes as before but with increasing frequency and mesh size inversely proportional to it. For OSM, we use 128, 256, or 512 subdomains on the smallest problem, but only 256 and 512 with ORAS, as 128 subdomains cannot fit in memory. As before, 64 sources are simulated, with full batch in case of OSM, and a batch size of 16 for ORAS. Each subdomain gets assigned 2~threads. \Cref{tab:time_osm} shows the result with OSM and \cref{tab:time_oras} with ORAS. In both tables we report time spent in local solves, in the GMRES orthogonalization process, and in sparse matrix products (the batched form of either the global matrix-vector product in ORAS, or the local SpMV in OSM). We also report the total time, which includes other minor components such as communications. In each case, the reported time is the maximum over all processes.

These tables show that OSM is consistently faster than ORAS, with a factor around 2. The absence of overlap allows for faster subdomain solves, the substructuring reduces the cost of GMRES, and replacing the global sparse matrix-vector product by local, interface-located operations all contribute to the speedup. The larger batch size also helps, although it could have been higher with ORAS with sufficient memory.

Despite the number of cores per unknown being kept stable, the scalability is not ideal, irrespective of the method. This well-known issue is due to the local nature of communications in one-level domain decomposition methods, and should be mitigated using two-level methods, which are out of the scope of this paper. For a fixed subdomain size, computation time seems to scale linearly with the frequency (especially with OSM), i.e., as the cubic root of the number of unknowns. Overall, this yields a complexity of $\mathcal{O}(f^{4})$ for the total amount of work. 

\begin{table}[ht]
\centering
\begin{tabular}{|c|c|c|c|c|c|c|c|c|}
\hline
\multicolumn{3}{|c|}{Problem solved} & \multicolumn{4}{c|}{Time spent in... (s)} & \multicolumn{2}{c|}{Cumulated} \\ \hline
$f$ (Hz) & MDOFs & $N$ & Solves & Ortho. & SpMM & Total & CPU-h & Its \\ \hline1.6 & 11.0 & 128 & 55.5 & 3.7 & 4.3 & 66.5 & 4.7 & 54 \\ \hline
1.6 & 11.0 & 256 & 31.2 & 3.3 & 2.6 & 41.1 & 5.8 & 70 \\ \hline
1.6 & 11.0 & 512 & 18.3 & 3.0 & 2.0 & 26.1 & 7.4 & 92 \\ \hline
2.0 & 21.1 & 256 & 69.1 & 5.0 & 5.8 & 85.3 & 12.1 & 72 \\ \hline
2.0 & 21.1 & 512 & 39.3 & 4.8 & 3.7 & 53.0 & 15.1 & 93 \\ \hline
2.0 & 21.1 & 1024 & 21.9 & 4.1 & 2.8 & 32.6 & 18.6 & 116 \\ \hline
2.5 & 41.0 & 512 & 97.8 & 8.5 & 8.3 & 122.0 & 34.7 & 104 \\ \hline
2.5 & 41.0 & 1024 & 53.2 & 6.7 & 5.3 & 72.7 & 41.4 & 128 \\ \hline
2.5 & 41.0 & 2048 & 30.7 & 6.0 & 3.8 & 45.9 & 52.3 & 166 \\ \hline
3.125 & 79.3 & 1024 & 123.3 & 11.4 & 11.2 & 155.1 & 88.2 & 136 \\ \hline
3.125 & 79.3 & 2048 & 77.1 & 9.7 & 7.3 & 104.9 & 119.4 & 177 \\ \hline
3.125 & 79.3 & 4096 & 45.0 & 9.8 & 5.9 & 68.8 & 156.6 & 252 \\ \hline
4.0 & 165.2 & 2048 & 170.4 & 16.4 & 15.9 & 217.6 & 247.6 & 179 \\ \hline
4.0 & 165.2 & 4096 & 105.0 & 14.8 & 10.5 & 143.6 & 326.7 & 245 \\ \hline
4.0 & 165.2 & 8192 & 59.1 & 12.2 & 7.6 & 89.1 & 405.5 & 312 \\ \hline
5.0 & 326.4 & 4096 & 226.3 & 22.9 & 21.6 & 291.1 & 662.5 & 249 \\ \hline
5.0 & 326.4 & 8192 & 130.5 & 18.8 & 14.7 & 182.4 & 830.2 & 318 \\ \hline
5.0 & 326.4 & 16384 & 70.5 & 13.3 & 8.6 & 99.5 & 905.6 & 346 \\ \hline
6.25 & 630.1 & 8192 & 303.1 & 29.4 & 28.7 & 378.6 & 1723.3 & 318 \\ \hline
6.25 & 630.1 & 16384 & 142.0 & 20.7 & 16.4 & 205.5 & 1870.7 & 357 \\ \hline
\end{tabular}
\caption{Results with OSM. MDOFs indiciates the number of millions of unknowns and $N$ the number of subdomains.}
\label{tab:time_osm}\end{table}
 \begin{table}[ht]
\centering
\begin{tabular}{|c|c|c|c|c|c|c|c|c|}
\hline
\multicolumn{3}{|c|}{Problem solved} & \multicolumn{4}{c|}{Time spent in... (s)} & \multicolumn{2}{c|}{Cumulated} \\ \hline
$f$ (Hz) & MDOFs & $N$ & Solves & Ortho. & SpMM & Total & CPU-h & Its \\ \hline
1.6 & 11.0 & 256 & 50.5 & 24.3 & 23.9 & 86.9 & 12.4 & 52 \\ \hline
1.6 & 11.0 & 512 & 31.5 & 13.2 & 16.0 & 50.5 & 14.4 & 64 \\ \hline
2.0 & 21.1 & 512 & 68.2 & 29.0 & 32.7 & 109.6 & 31.2 & 71 \\ \hline
2.0 & 21.1 & 1024 & 45.0 & 22.6 & 21.4 & 71.7 & 40.8 & 89 \\ \hline
2.5 & 41.0 & 1024 & 91.0 & 42.8 & 46.0 & 149.0 & 84.8 & 97 \\ \hline
2.5 & 41.0 & 2048 & 55.0 & 25.5 & 27.3 & 87.9 & 100.0 & 116 \\ \hline
3.125 & 79.3 & 2048 & 129.3 & 64.9 & 60.5 & 197.5 & 224.7 & 128 \\ \hline
3.125 & 79.3 & 4096 & 75.4 & 37.2 & 35.9 & 120.2 & 273.5 & 161 \\ \hline
4.0 & 165.2 & 4096 & 173.4 & 73.4 & 85.5 & 277.7 & 631.9 & 177 \\ \hline
4.0 & 165.2 & 8192 & 111.1 & 55.8 & 53.8 & 176.1 & 801.4 & 225 \\ \hline
5.0 & 326.4 & 8192 & 223.8 & 110.6 & 98.5 & 363.7 & 1655.2 & 244 \\ \hline
5.0 & 326.4 & 16384 & 160.6 & 89.2 & 82.6 & 256.9 & 2338.2 & 344 \\ \hline
6.25 & 630.1 & 16384 & 311.1 & 162.3 & 147.2 & 502.6 & 4574.7 & 347 \\ \hline
\end{tabular}
\caption{Results with ORAS. MDOFs indiciates the number of millions of unknowns and $N$ the number of subdomains.}
\label{tab:time_oras}\end{table}
 
\pgfplotsset{compat=1.18}
\usepgfplotslibrary{groupplots}
\begin{figure}[ht]\centering
\begin{tikzpicture}
\begin{groupplot}[
  group style={group size=2 by 1, horizontal sep=2cm},
  width=6cm, height=6cm,
  ymin=0, ymax=600,
  ylabel={Solving time (s)},
  xlabel={Frequency (Hz)},
  grid=both,
  legend to name=sharedlegend,
  legend columns=1,
  legend style={/tikz/every even column/.append style={column sep=1cm}},
  legend image post style={mark options={solid}},
  every mark/.append style={solid},
]
\nextgroupplot[title={ORAS}]
\addplot+[thick, mark=*, color=blue] coordinates {};
\addplot+[thick, mark=square*, color=blue] coordinates {
  (1.6,86.89)   (2,109.6)   (2.5,149)   (3.125,197.5)   (4,277.7)   (5,363.7)   (6.25,502.6) };
\addplot+[thick, mark=triangle*, color=red] coordinates {
  (1.6,50.47)   (2,71.65)   (2.5,87.89)   (3.125,120.2)   (4,176.1)   (5,256.9) };
\addplot[dashed, domain=1.5:6.5] {{89.07986989419483*x + -70.4597450943455}};
\addplot[dashed, domain=1.5:6.5] {{59.63133371363376*x + -53.93870948849579}};
\nextgroupplot[title={OSM}]
\addplot+[thick, mark=*, color=green!60!black] coordinates {
  (1.6,66.46)   (2,85.27)   (2.5,122)   (3.125,155.1)   (4,217.6)   (5,291.1)   (6.25,378.6) };
\addplot+[thick, mark=square*, color=blue] coordinates {
  (1.6,41.11)   (2,52.97)   (2.5,72.74)   (3.125,104.9)   (4,143.6)   (5,182.4)   (6.25,205.5) };
\addplot+[thick, mark=triangle*, color=red] coordinates {
  (1.6,26.06)   (2,32.64)   (2.5,45.94)   (3.125,68.83)   (4,89.1)   (5,99.49) };
\addplot[dashed, domain=1.5:6.5] {{67.81542519858088*x + -49.07090453360957}};
\addlegendentry{Large subdomains}
\addplot[dashed, domain=1.5:6.5] {{37.72987925025901*x + -17.174742092869906}};
\addlegendentry{Medium subdomains}
\addplot[dashed, domain=1.5:6.5] {{23.24463266780757*x + -10.262438395132136}};
\addlegendentry{Small subdomains}
\end{groupplot}
\node at ($(group c1r1.south)!0.5!(group c2r1.south)-(0,2.5cm)$) {\pgfplotslegendfromname{sharedlegend}};
\end{tikzpicture}
\caption{Weak scaling comparison of OSM and ORAS with constant subdomain size. Dashed lines represent affine functions fitting each curve.}
\end{figure}
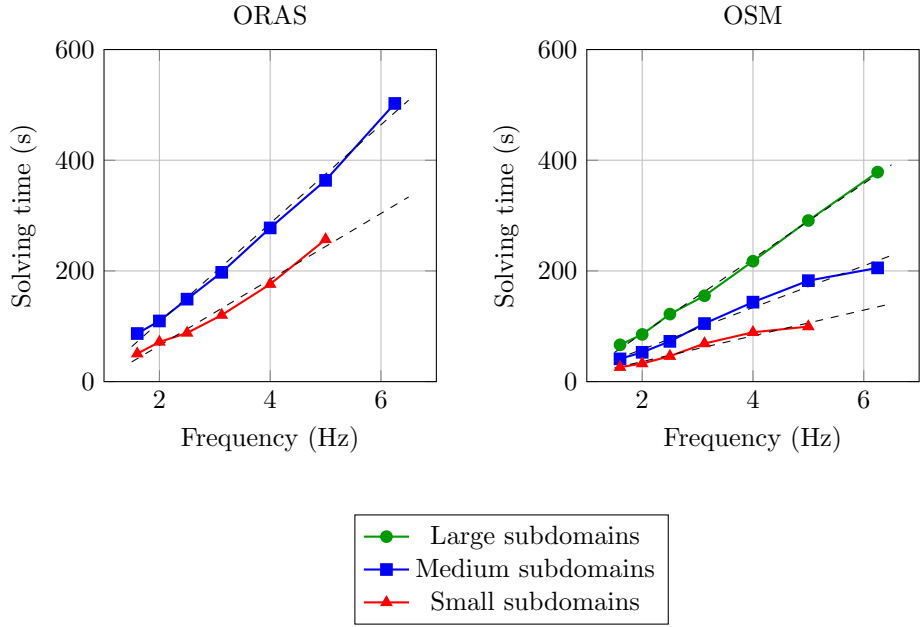
 
\paragraph{Comparison with a direct solver}

We also solved the system using MUMPS as a direct solver with BLR enabled ($\varepsilon_{BLR} = 10^{-4}$). For the \SI{1.6}{\hertz} case (about 11 million DOFs) with 256 domains (and 512 threads, on 4 nodes), a global factorization is computed in \SI{74}{\s}, and then \SI{6}{s} are needed to solve the 64 sources, leading to a relative residual around $10^{-3}$. With the same resources, OSM achieved a solution in 41 seconds. A direct solve is thus slower but has a smaller marginal cost per RHS, which could make it competitive when the number of sources is huge. A simple affine extrapolation shows that for this small-size problem the break-even point is around 150 sources, which is not rare in real FWI cases, especially when the adjoint problems are taken into account, or when the factorization is often reused, such as with second-order adjoint methods~\cite{metivier:hal-00826614} used in FWI. However, it is known that standard multifrontal direct solvers such as MUMPS do not scale as well as domain decomposition methods with problem size: full-rank MUMPS has a quadratic complexity with respect to the number of unknowns for 3D problems~\cite{mary2017}, which is however partially mitigated by the BLR approach.
In practical FWI applications, using direct solves on low frequencies and iterative methods on high frequencies thus seems to be the best compromise, as already suggested in~\cite{operto_is_2023}.

\begin{figure}
\centering
\begin{tikzpicture}
  \begin{axis}[
    ybar stacked,
    bar width=22pt,
    height=7cm,
    width=6cm,
    ylabel={Time / Cost},
    symbolic x coords={ORAS,OSM},
    xtick=data,
    enlarge x limits=0.25,
    legend style={at={(0.5,-0.18)},anchor=north,legend columns=-1},
    every axis plot/.style={/tikz/line width=0pt},  ]
\addplot[fill=blue] coordinates
      {(ORAS,162.3) (OSM,20.7)};

\addplot[fill=red] coordinates
      {(ORAS,147.2) (OSM,16.4)};

\addplot[fill={green!60!black}] coordinates
      {(ORAS,311.1) (OSM,142.0)};

    \legend{GMRES orthogonalization,SpMM,Local solves}
  \end{axis}
\end{tikzpicture}
    \caption{Time repartition in the OSM (left) and ORAS (right) methods for the 600M unknowns problem on \num{16384} domains. Actual total times are \SI{502.6}{\s} and \SI{205.5}{\s} respectively. }
\label{fig:stack_chart}
\end{figure}
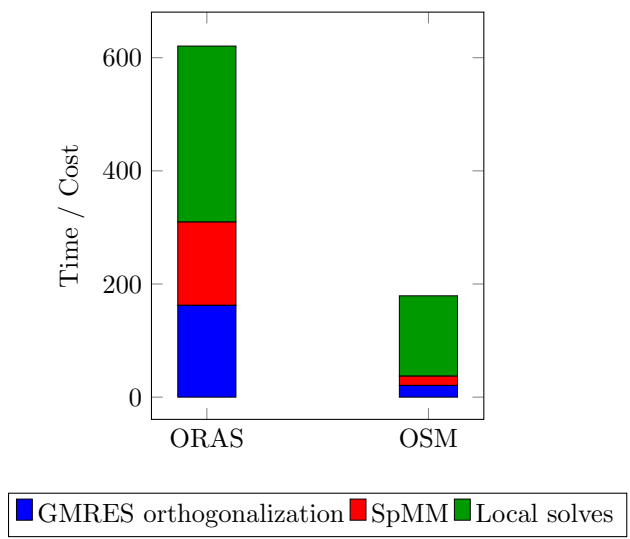

\paragraph{Computation time breakdown during the solving phase}
For both methods, the main component of the computation (once the subproblem matrices have been factorized) is the repeated solution of each subproblem. However, a significant amount of time is also spent on two additional tasks. First, GMRES needs to build an orthonormal basis of the Krylov subspace using a Gram--Schmidt procedure. Second, at each iteration, the right-hand side of the subproblems must be computed. In ORAS, they are the restriction of the vector to which the preconditioner is applied, i.e., the product of the FEM matrix $A$ and the latest Krylov vector. In OSM, it is the inhomogeneous source coming from other subdomains which must be assembled. In both cases, this step can be represented by sparse matrix-vector products. In OSM, these matrices are smaller (and rectangular) since they map interface variables to the source of the local Helmholtz problem. Thus, there are three main components to consider and all of which have smaller dimensions in OSM, which explains the superior performance of the method. A more detailed analysis presented in \cref{fig:stack_chart} shows which portion of the time is spent in these three parts for our largest test case. As expected, the GMRES part and the sparse algebra are radically cheaper in OSM due to the substructuring, allowing the algorithm to spend most of its resources on the local solves. Combined with the smaller subdomains that are faster to solve, the whole method becomes significantly more efficient. 

\section{Conclusion and perspectives}
\label{sec:conclusion}

We tested and compared optimized implementations of the ORAS preconditioner and the OSM non-overlapping DDM for the Helmholtz equation, with a focus on solving large-scale geophysical problems with multiple right-hand sides, using one of the world's largest supercomputer, LUMI. We first showed that, while ORAS has a substantially faster convergence with basic impedance conditions, a smarter choice of parameters allows OSM to almost bridge the gap. With optimized transmission conditions the OSM method becomes consistently faster and more memory-efficient than ORAS: having smaller subdomain sizes, applying GMRES to the substructured equation, and replacing the global sparse matrix-vector products by localized interface operations all contribute to this difference. Furthermore, OSM reduces the memory footprint of large batches, which can provide additional speedup.

While OSM is quite clearly the most appropriate one-level method in this case, two-level methods must be considered for larger problems, to solve the scalability issues. Geometric coarse-grids were successfully used for ORAS~\cite{tournier:hal-03942570} and have a moderate cost (essentially limited to storing the coarse problems), and recent progress has been made on spectral coarse grids~\cite{dolean2024schwarzpreconditionerhkgeneocoarse, coarse_bootland, ConenDTN}, which however are designed for overlapping methods. Developing a two-level method for OSM (either geometric or spectral) is a promising avenue for future work and is currently being investigated.

Finally, similar investigations should be carried for other equations, such as Maxwell or elasticity, where the same principles apply. OSM should be more economical than ORAS per-iteration for the same reasons (smaller domains, substructuring, no global sparse matrix-vector product), but a convergence analysis is needed to ensure the difference in iteration count does not outweigh these advantages.

\section{Reproducibility}
\label{sec:reproducibility}

The source code to reproduce the numerical results is available at \url{https://gitlab.onelab.info/boris-martin/oras_vs_osm}. Instructions for building and running the tests are available in the \href{https://gitlab.onelab.info/boris-martin/oras_vs_osm/-/blob/master/README.md}{README.md} file.

\section*{Acknowledgments}

The first author gratefully acknowledges financial support from the Fonds de la Recherche Scientifique – FNRS (Belgium) through a FRIA doctoral fellowship. The research was funded in part by the Walloon Region through the Win2Wal EXPANSION project (Grant No. 2010161). Computational resources have been provided by the Consortium des Équipements de Calcul Intensif (CÉCI), funded by the Fonds de la Recherche Scientifique de Belgique (F.R.S.-FNRS) under Grant No. 2.5020.11 and by the Walloon Region.
The present research benefited from computational resources made available on Lucia, the Tier-1 supercomputer of the Walloon Region, infrastructure funded by the Walloon Region under the grant agreement n\textdegree 1910247.
We acknowledge LUMI-BE\footnote{LUMI-BE is joint effort from BELSPO (federal), SPW Économie, Emploi, Recherche (Wallonia), Department of Economy, Science \& Innovation (Flanders) and Innoviris (Brussels).} for awarding this project access to the LUMI supercomputer, owned by the EuroHPC Joint Undertaking, hosted by CSC (Finland) and the LUMI consortium through a LUMI-BE Regular Access call.

\clearpage
\bibliographystyle{siamplain}
\bibliography{references, BLR}
\end{document}